\documentclass[12pt,a4paper]{article}
\usepackage[T1]{fontenc}
\usepackage[scaled=0.91]{helvet}
\usepackage{courier}
\usepackage{amsmath,amssymb,amsfonts}        % moduli per le estensioni ed i fonts AMS
\usepackage{mathptmx}
\usepackage[dvips]{graphicx}
\usepackage{psfrag}
\usepackage{color}
\usepackage[pdfauthor={Riccardo Fazio},pdftitle={American Put Option}, bookmarks,colorlinks]{hyperref}

\def\ds{\displaystyle}
\def\FBP{free boundary problem}
\def\RR{\hbox{I\kern-.2em\hbox{R}}}
\def\NN{\hbox{I\kern-.2em\hbox{N}}}
\def\ds{\displaystyle}

\newcommand{\eqnsection}{
   \renewcommand{\theequation}{{\thesection.\arabic{equation}}}
   \makeatletter
   \csname @addtoreset\endcsname{equation}{section}
   \makeatother}
\eqnsection

\title{Perpetual American Put Option: \\
an Error Estimator for a\\
Non-Standard Finite Difference Scheme}
\author{Riccardo Fazio\thanks{Author home-page and e-mail: http://mat521.unime.it/$\sim$fazio \ \ \ rfazio@unime.it}\\
Department of Mathematics and Computer Science \\
University of Messina \\
Viale F. Stagno D'Alcontres 31, 98166 Messina, Italy}
\date{\today}
\begin{document}
\maketitle
\begin{abstract}
In this paper we present a MATLAB version of a non-standard finite difference scheme for the numerical solution of the perpetual American put option models of financial markets. 
These models can be derived from the celebrated Black-Scholes models letting the time goes to infinity.
The considered problem is a \FBP \ defined on a semi-infinite interval, so that it is a non-linear problem complicated by a boundary condition at infinity.
By using non-uniform maps, we show how it is possible to apply the boundary condition at infinity exactly.
Moreover, we define a posteriori error estimator that is based on Richardson's classical extrapolation theory.
Our finite difference scheme and error estimator are favourably tested for a simple problem with a known exact analytical solution.
\end{abstract}
\smallskip

\noindent
{\bf Key Words:} 
Perpetual American put option, free boundary problem, non-standard finite difference scheme, posteriori error estimator. 
%\smallskip

%\noindent
%{\bf AMS Subject Classifications:} .

%\newpage

%\textcolor{red}{Bigliograpy}
\section{Introduction}
Analytical solutions of models of American option problems are seldom available, so such derivatives of financial markets must be priced by numerical methods (Amin and Khanna \cite{Amin:1994:CAO}, Barraquand and Pudet \cite{Barraquand:1994:PAP}, Broadie and Detemple \cite{Broadie:1996:MVR}, Nielsen et al. \cite{Nielsen:2002:PFF}, Barone-Adesi \cite{Baroni:2005:SAP}, D\"{u}ring and Fourni\'e \cite{During:2012:HOC} or Milev and Tagliani \cite{Nilev:2013:EIS}).
In this paper we present a MATLAB version of a non-standard finite difference scheme for the numerical solution of the perpetual American put option models of financial markets. 
These models can be derived from the celebrated Black-Scholes models (Leland \cite{Leland:1985:OPR}, Avellaneda and Par\'as \cite{Avellaneda:1996:MVR},Frey and Patie \cite{Frey:2002:RMD} and Janda\v{c}ka and \v{S}ev\v{c}ovi\v{c} \cite{Jandacka:2005:RAP})
 letting the time goes to infinity (Bensoussan \cite{Bensoussan:1984:TOP} or Elliot and Kopp \cite[pp. 196-199]{Elliot:1999:MFM}).
The considered problem is a \FBP \ defined on a semi-infinite interval, so that it is a non-linear problem complicated by a boundary condition at infinity.
By using non-uniform maps, we show how it is possible to apply the boundary condition at infinity exactly.
Non-uniform maps have been applied to the numerical solution of ordinary and partial differential equations on unbounded domains (van de Vooren and Dijkstra \cite{vandeVooren:1970:NSS}, Botta et al. \cite{Botta:1972:NSN}, Davis \cite{Davis:1972:NSN}, Grosch and Orszag \cite{Grosch:NSP:1977}, Boyd \cite{Boyd:2001:CFS}, Koleva \cite{Koleva:NSH:2006} or Fazio and Jannelli \cite{Fazio:2014:FDS}).
Moreover, we deduce a posteriori error estimator within Richardson's classical extrapolation theory.
Our finite difference scheme and error estimator are favourably tested for a simple problem with a known exact analytical solution.
From the obtained numerical results we can asses that: the finite difference method is second order accurate, the numerical solution can be improved by repeated Richardson's extrapolations and the error estimator provides upper bounds for the exact error.

\section{Perpetual American put option}\label{S:PPO}
In order to test our error estimator, in this section, we consider a test problem with known exact analytical solution.
This problem is a \FBP \ arising as a simple toy model in the study of financial markets \cite{Bensoussan:1984:TOP}.
A mathematical model describing the perpetual American put option is given by
\begin{align}\label{PPO:model}
& \frac{1}{2} \sigma^2 S^2  
{\ds \frac{d^2 P}{dS^2}} + r S {\ds \frac{dP}{dS}} -r P = 0 \ , \qquad \mbox{on} \ \ R \le S < \infty \ , \nonumber \\
& P(R) = \max\{E - R, 0\} \ , \qquad {\ds \frac{dP}{dS}}(R) = - 1 \ , \\ %&  P(S, T) = \max(E-S, 0) \ , \qquad R(T) = E \ ,  \nonumber \\
& \lim_{S\rightarrow \infty} P(S) = 0 \ , \nonumber 
\end{align}
where $S$ is the price of a given asset, $P(S)$ is the price of the perpetual American put option to sell the asset, $R$ is the unknown free boundary, $\sigma$, $r$ and $E$ are the volatility, interest rate and exercise price of the asset, respectively.
This problem (\ref{PPO:model}) has the exact solution
\begin{equation}
P(S) = (E-R)\; R^{2\; r/\sigma^2}\; S^{-2\; r/\sigma^2} \ , \qquad R = {\ds \frac{2\; r\; E}{2\; r + \sigma^2}} \ ,
\end{equation}
see \cite[pp. 196-199]{Elliot:1999:MFM}.
In order to fix the domain, see Crank \cite[pp. 187-192]{Crank:FMB:1984}, we can apply Landau's transformation of variables
\[
x = S/R \ , \qquad u(x) = P(x \; R) \ .
\]
In the new variables the put option problem (\ref{PPO:model}) can be rewritten as follows
\begin{align}\label{APO:model:fix}
& \frac{1}{2} \sigma^2 x^2  
{\ds \frac{d^2 u}{dx^2}} + r x {\ds \frac{du}{dx}} -r u = 0 \ , \qquad \mbox{on} \ \ 1 \le x < \infty \ , \nonumber \\
& u(1) = \max\{E - R, 0\} \ , \qquad {\ds \frac{du}{dx}}(1) = - R \ , \\ %&  P(S, T) = \max(E-S, 0) \ , \qquad R(T) = E \ ,  \nonumber \\
& \lim_{x\rightarrow \infty} u(x) = 0 \ , \nonumber 
\end{align}
Moreover, this model can be rewritten in standard form as a first order system of ordinary differential equations.
The model (\ref{PPO:model}) is a special instance of the American put option obtained formally by letting the time variable to go to infinity.
In recent years several generalization, ranging from the introduction of further relevant markets parameters to non-constants volatility and the like, of this model have been proposed in literature.
In particular, one can take into account: the presence of transaction costs (see e.g. Leland \cite{Leland:1985:OPR}, Avellaneda and Par\'as \cite{Avellaneda:1996:MVR}), feedback and illiquid market effects due to large traders choosing given stock-trading strategies (Frey and Patie \cite{Frey:2002:RMD}), risk from unprotected portfolio (Janda\v{c}ka and \v{S}ev\v{c}ovi\v{c} \cite{Jandacka:2005:RAP}).
In order to take into account also those different models, and using the fixed boundary formulation (\ref{APO:model:fix}), we study here the following class of problems 
\begin{align}\label{S:model}
& {\ds \frac{du}{dx}} = f(x,u,v) \ , \qquad  \mbox{on} \ \ 1 \le x < \infty \ , \nonumber \\[1.2ex]
& {\ds \frac{dv}{dx}} = g(x,u,v) \ ,  \nonumber \\[-1ex]
& \\[-1ex]
& {\ds \frac{dR}{dx}} = 0 \ ,  \nonumber \\[1ex] 
& u(1) = \max\{E - R, 0\} \ , \qquad v(1) = - R \ ,  %&  P(S, T) = \max(E-S, 0) \ , \qquad R(T) = E \ ,  \nonumber \\
\qquad \lim_{x\rightarrow \infty} u(x) = 0 \ , \nonumber 
\end{align}
where $R$ is treated as a supplentary variable because its value is unknown and has to be found as part of the solution.
Of coarse, our benchmark problem (\ref{APO:model:fix}) belongs to (\ref{S:model}) for a suitable change of variables and suitable functional form of $f$ and $g$.

\section{Quasi-uniform grids}\label{S:quniform}
Let us consider the smooth strict monotone quasi-uniform maps $x = x(\xi)$, the so-called grid generating functions,
\begin{equation}\label{eq:qu1}
x = -c \cdot \ln (1-\xi) + 1\ ,
\end{equation}
and
\begin{equation}\label{eq:qu2}
x = c \frac{\xi}{1-\xi} + 1\ ,
\end{equation}
where $ \xi \in \left[0, 1\right] $, $ x \in \left[1, \infty\right] $, and $ c > 0 $ is a control parameter.
So that, a family of uniform grids $\xi_n = n/N$ defined on interval $[0, 1]$ generates one parameter family of quasi-uniform grids $x_n = x (\xi_n)$ on the interval $[1, \infty]$.
The two maps (\ref{eq:qu1}) and (\ref{eq:qu2}) are referred as logarithmic and algebraic map, respectively. 
The logarithmic map (\ref{eq:qu1}) gives slightly better resolution near $x = 1$ than the
algebraic map (\ref{eq:qu2}), while the algebraic map gives much better resolution than the
logarithmic map as $x \rightarrow \infty$. 
In fact, it is easily verified that
\[
-c \cdot \ln (1-\xi) + 1 < c \frac{\xi}{1-\xi} +1 \ ,
\]
for all $\xi$.

The problem under consideration can be discretized by introducing a uniform grid $ \xi_n $ of $N+1$ nodes in $ \left[0, 1\right] $ with $\xi_0 = 0$ and $ \xi_{n+1} = \xi_n + h $ with $ h = 1/N $, so that $ x_n $ is a quasi-uniform grid in $ \left[1, \infty\right] $. 
The last interval in (\ref{eq:qu1}) and (\ref{eq:qu2}), 
namely $ \left[x_{N-1}, x_N\right] $, is infinite but the point $ x_{N-1/2} $ is finite, because the non integer nodes are defined by 
\[
x_{n+\alpha} = x\left(\xi=\frac{n+\alpha}{N}\right) \ ,
\]
with $ n \in \{0, 1, \dots, N-1\} $ and $ 0 < \alpha < 1 $.
%In this way we have defined also the ghost cell $ x_{-1} $.  
These maps allow us to describe the infinite domain by a finite number of intervals.
The last node of such grid is placed on infinity so right boundary conditions
are taken into account correctly.

\section{A non-standard finite difference scheme}\label{S:NSFDS}
We can approximate the values of $u(x)$ on the mid-points of the grid
\begin{equation}
u_{n+1/2} \approx \frac{x_{n+3/4}-x_{n+1/2}}{x_{n+3/4}-x_{n+1/4}} u_n + \frac{x_{n+1/2}-x_{n+1/4}}{x_{n+3/4}-x_{n+1/4}} u_{n+1} \ .
\label{eq:u:mod}
\end{equation}
that is, a non-standard central difference formula.
Taking into account the results by Veldam and Rinzema \cite{Veldam:PNG:1992}, for the first derivative at the mid-points of the grid we can apply the following approximation
\begin{equation}
\frac{du}{dx}(x_{n+1/2}) \approx \frac{u_{n+1}-u_n}{2\left(x_{n+3/4} - x_{n+1/4}\right)} \ ,
\label{eq:du}
\end{equation}
that is, again, a non-standard central difference formula.
These finite difference formulae use the value $ u_N = u_\infty $, but not $ x_N = \infty $.
The approximation (\ref{eq:u:mod}) is a variant of the formula used by Fazio and Jannelli \cite{Fazio:2014:FDS}. 
A non-standard finite difference scheme on a quasi-uniform grid for our financial problem
(\ref{PPO:model}) can be defined by using the approximations given by (\ref{eq:u:mod}) and (\ref{eq:du}) above.

We denote by the $3-$dimensional  vector $ {\bf U}_n = (U_n, V_n, R)^T $ the numerical approximation to the solution $ {\bf u} (x_n) = (u(x_n),v(x_n),R)^T $ of (\ref{S:model}) at the points of the mesh, that is for $ n = 0, 1, \dots , N $ .
A finite difference scheme for
(\ref{S:model}) can be written as follows:
\begin{eqnarray}\label{boxs}
& {\bf U}_{n+1} - {\bf U}_{n} - a_{n+1/2} {\bf f} \left( x_{n+1/2}, b_{n+1/2}{\bf U}_{n+1} + c_{n+1/2}{\bf U}_{n} \right) = {\bf 0}
\ , \nonumber\\[-1.3ex]
& \\[-1.3ex]
& {}_1 {\bf U}_{0} = \max\{E-R, 0 \} \ , \quad {}_2 {\bf U}_{0} = -R \ , \quad {}_1 {\bf U}_{N}  = 0 \ ,  \nonumber
\end{eqnarray}
for $n=0$, $1$, $\dots$ , $N-1$, here ${\bf f} = (f, g, 0)^T$, ${}_j {\bf U}$ is the $j$-component of the vector ${\bf U}$ and
\begin{eqnarray}\label{eq:abc}
a_{n+1/2} &=& 2\left(x_{n+3/4} - x_{n+1/4}\right) \ , \nonumber \\
b_{n+1/2} &=& \frac{x_{n+1/2}-x_{n+1/4}}{x_{n+3/4}-x_{n+1/4}} \ , \\
c_{n+1/2} &=& \frac{x_{n+3/4}-x_{n+1/2}}{x_{n+3/4}-x_{n+1/4}} \nonumber \ , 
\end{eqnarray}
for $n=0, 1, \dots , N-1$.

It is evident that (\ref{boxs}) is a nonlinear system of $ 3 \cdot (N+1)$ equations in the $ 3 \cdot (N+1)$ unknowns $ {\bf U} = ({\bf U}_0,
{\bf U}_1, \dots , {\bf U}_N)^T $. 

\section{Richardson's extrapolation}\label{S:extra}
The utilization of a quasi-uniform grid allows us to improve our numerical results.
The algorithm is based on Richardson's extrapolation, introduced by Richardson in \cite{Richardson:1910:DAL,Richardson:1927:DAL}, and it is the same for many finite difference methods: for numerical differentiation or integration, solving systems of ordinary or partial differential equations.
To apply Richardson's extrapolation, we carry on several calculations on embedded uniform or quasi-uniform grids with total number of nodes $N$: e.g., for the numerical results reported in the next section we used $N = 2$, $4$, $8$, $16$, $32$, $64$, $128$, $256$, $512$ or $N = 5$, $10$, $20$, $40$, $80$, $160$, $320$, $640$, $1280$. 
We can identify these grids with the index $g=0$, the coarsest one, $1$, $2$, and so on towards the finest grid.
Between two adjacent grids all nodes of largest steps are identical to even nodes of denser grid due to quasi-uniformity. 
To find an approximation of a scalar value $U$ we can apply $k$ Richardson's extrapolations on the used grids
\begin{equation}\label{eq:Rextra}
U_{g+1,k+1} = U_{g+1,k} + \frac{U_{g+1,k}-U_{g,k}}{q^{p_k}-1} \ ,
\end{equation}
where $g \in \{0, 1, 2 , \dots , G-1\}$, $k \in \{0, 1, 2, \dots , G-1\}$, $q = N_g/N_{g-1}$ is the grid refinement ratio, and $p_k$ is the true order of the discretization error, see Schneider and Marchi \cite{Schneider:GRR:2005} and the references quoted therein.
This formula is asymptotically exact in the limit as $N$ goes to infinity if we use uniform or quasi-uniform grids.
We notice that to obtain each value of $U_{g+1,k+1}$ requires having computed two solution $U$ in two adjacent grids, namely $g+1$ and $g$ at the extrapolation level $k$.
Hence, it gives the real value of numerical solution error without knowledge of exact solution.
For any $g$, the level $k=0$ represents the numerical solution of $U$ without any extrapolation, which is obtained as described in section \ref{S:NSFDS}.
The case $k=1$ is the classical single Richardson's extrapolation, which is usually used to estimate the discretization error or to improve the solution accuracy.
If we have computed the numerical solution on $G+1$ nested grids then we can apply equation (\ref{eq:Rextra}) $G$ times performing $G$ Richardson's extrapolation. 

Here we are interested to show how within Richardson's extrapolation theory we can derive an error estimate.
For any value of interest $U$, the numerical error $E$ can be defined by
\begin{equation}\label{eq:GE}
E = u - U \ ,
\end{equation}
where $u$ is the exact analytical solution.
Usually, we have several different sources of errors: discretization, round-off, iteration and programming errors.
Discretization errors are due to our replacement of a continuous problem with a discrete one and is errors can be reduced by reducing the discretization parameters, enlarging the value of $N$ in our case.
Round-off error are due to the utilization of floating-point arithmetic to implement the algorithms available to solve the discrete problem.
This kind of error can be reduced by using higher precision arithmetic, double or, when available, fourth precision.
Iteration errors are due to stopping an iteration algorithm that is converging but only as the number of iterations goes to infinity.
Of course, we can reduce this kind of error by requiring more restrictive termination criteria for our iterations, the iterations of \texttt{fsolve} MATLAB routine in the present case.    
Programming errors are behind the scope of this work but they can be eliminated or at least reduced by adopting what is called structured programming.  
When the numerical error is caused prevalently by the discretization error and in the case of smooth enough solutions the discretization error can be decomposed into a sum of powers of the inverse of $N$ 
\begin{equation}\label{eq:asymE}
u = U_N + C_0 \left(\frac{1}{N}\right)^{p_0}+ C_1 \left(\frac{1}{N}\right)^{p_1}+ C_2 \left(\frac{1}{N}\right)^{p_2}+ \cdots \ ,
\end{equation}
where $C_0$, $C_1$, $C_2$, $\dots$ are coefficient that depend on $u$ and its derivatives, but are independent on $N$, and $p_0$, $p_1$, $p_2$, $\dots$ are the true orders of the error.
The value of each $p_k$ is usually a positive integer with $p_0 < p_1 < p_2 < \cdots$ and constitute an arithmetic progression of ratio $p_1-p_0$.
The value of $p_0$ is called the asymptotic order or the order of accuracy of the method or of the numerical solution $U$. 
So that, the theoretical order of accuracy of the numerical solution $U$ with $k$ extrapolations the $p_k$ orders verify the relation
\begin{equation}\label{eq:pk}
p_k = p_0 + k (p_1-p_0) \ ,
\end{equation}
where this equation is valid for $k \in \{0, 1, 2, \dots , G-1\}$.
 
\subsection{Error estimate}\label{SS:error}
To show how Richardson's extrapolation can be also used to get an error estimate for the computed numerical solution we use the notation introduced above.
By replacing into equation (\ref{eq:asymE}) $N$ with $2N$ and subtracting, to the obtained equation, equation (\ref{eq:asymE}) times $(1/2)^{p_0}$ we get the first extrapolation formula 
\begin{equation}\label{eq:Rextra1}
u \approx  U_{2N} + \frac{U_{2N}-U_{N}}{2^{p_0}-1} \ ,
\end{equation}
that has a leading order of accuracy equal to $p_1$.
Taking into account equation (\ref{eq:Rextra1}) we can conclude that the error estimate by a first Richardson's extrapolation is given by
\begin{equation}\label{eq:est1}
Eest = \frac{U_{2N}-U_{N}}{2^{p_0}-1} \ ,
\end{equation}
where $p_0$ is the order of the numerical method used to compute the numerical solutions.
In comparison with (\ref{eq:est1}) a safer error estimator can be defined by
\begin{equation}\label{eq:est2}
Esafe = U_{2N}-U_{N} \ .
\end{equation}
Of course, $p_0$ can be found by
\begin{equation}\label{eq:p}
p_0 \approx {\ds \frac{\log(|U_N-u|)-\log(|U_{2N}-u|)}{\log(2)}} \ ,
\end{equation}
where $u$ is again the exact solution (or, if the exact solution is unknown, a reference solution computed with a suitable large value of $N$), and both $u$ and $U_{2N}$ are evaluated at the same grid-points of $U_N$. 

\section{Numerical results}
It should be mentioned that all numerical results reported in this paper were performed on an ASUS personal computer with i7 quad-core Intel processor and 16 GB of RAM memory running Windows 8.1 operating system. 

The non-standard finite difference scheme described above has been implemented in MATLAB.
In this way we take advantage of the available MATLAB built-in functions.
In particular, for the solution of the non-linear system (\ref{boxs}) we used the function \texttt{fsolve}.
Among the available alternative we used the \lq \lq Levenberg-Marquardt\rq \rq \ with $TolFun=10^{-15}$ and $TolX=10^{-15}$ options.
These values of $TolFun$ and $TolX$ define the termination criteria for \texttt{fsolve}. 
Usually, the \texttt{fsolve} routine took between 5 to 11 iterations to get a numerical solution that verifies the stopping criteria. 

\begin{figure}[!hbt]
\centering
\psfrag{x}[][]{\small $S$}
\psfrag{u}[][]{\small $$} 
\psfrag{u1}[][]{\small $P_n$} 
\psfrag{u2}[][]{\small $\frac{dP}{dS}|_n$} 
\psfrag{ue1}[][]{\small $P(S_n)$} 
\psfrag{ue2}[][]{\small $\frac{dP}{dS}({S_n})$} 
%\hfil \includegraphics[width=.4\textwidth]{PPO_exact} \hfil
%\includegraphics[width=.4\textwidth]{QU_PPO_c20_N160} \hfil
\includegraphics[width=.9\textwidth]{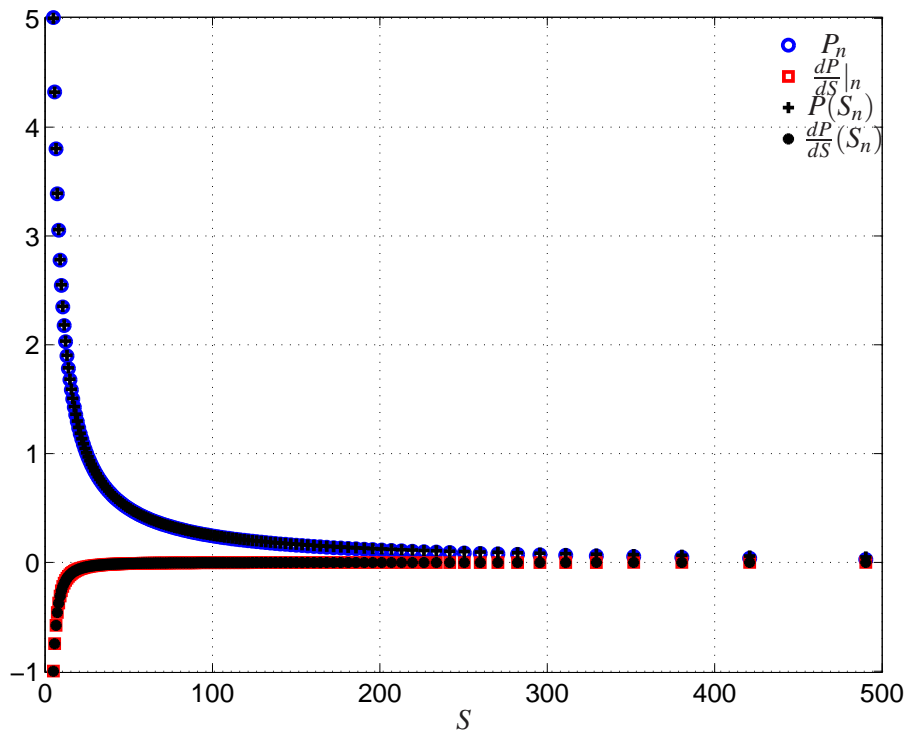}
\caption{\it Sample exact and numerical solutions for the test problem (\ref{PPO:model}).%}
The symbols indicate: $+$ the exact solution, $\bullet$ its first derivative, $\circ$ the numerical solution, and $\square$ its numerical first derivative.}
\label{fig:APO:exact}
\end{figure}

To set a specific test problem we fixed the following values for the involved parameters 
\begin{equation}\label{eq:parameters}
\sigma^2 = 0.1 \ , \qquad r = 0.05 \ , \qquad E = 10 \ .
\end{equation}
As we will see below these values provides an exact solutions that remains different from zero within a large domain. 
For our numerical computations we used both the two maps (\ref{eq:qu1}) with $c=20$ and (\ref{eq:qu2}) with $c=10$, but the results reported below are concerned with the fist map because the results obtained with the second map are, indeed, very similar.
In order to speed up the computations for different values of $N$ we adopted a continuation strategy.
For a small value of $N$, usually $N=2$ or $N=5$, we always used a constant initial iterate vector made with all components equal to one.
Then, when refining the grid we used the accepted final iterate of the previous value of $N$ as first iterate for the computation with the next value of $N$.  
Figure \ref{fig:APO:exact} shows a reference solution.
The numerical results for $N = 128$ can be seen on the same figure.
%We can note the uniform grid on the left frame and the non-uniform one on the right frame where the last grid-point (not shown) is located at infinity.  
The non-uniform grid is clearly visible even if the last grid-point is not shown because it is located at infinity.  

\begin{figure}[!hbt]
\centering
\psfrag{x}[][]{\small $S$} 
\psfrag{E}[][]{\small $E(S), Esafe$} 
\psfrag{E1}[][]{\small $E_1$} 
\psfrag{E2}[][]{\small $E_2$} 
\psfrag{Esafe1}[][]{\small $Esafe_1$} 
\psfrag{Esafe2}[][]{\small $Esafe_2$} 
%\includegraphics[width=.7\textwidth]{QU_APO_error_c20_N160} \\
%\hfil \includegraphics[width=.4\textwidth]{QU_PPO_E_c20_N20} \hfil
%\includegraphics[width=.4\textwidth]{QU_PPO_E_c20_N40} \hfil \\
%\hfil \includegraphics[width=.4\textwidth]{QU_PPO_E_c20_N80} \hfil
%\includegraphics[width=.4\textwidth]{QU_PPO_E_c20_N160} \hfil
\includegraphics[width=.7\textwidth]{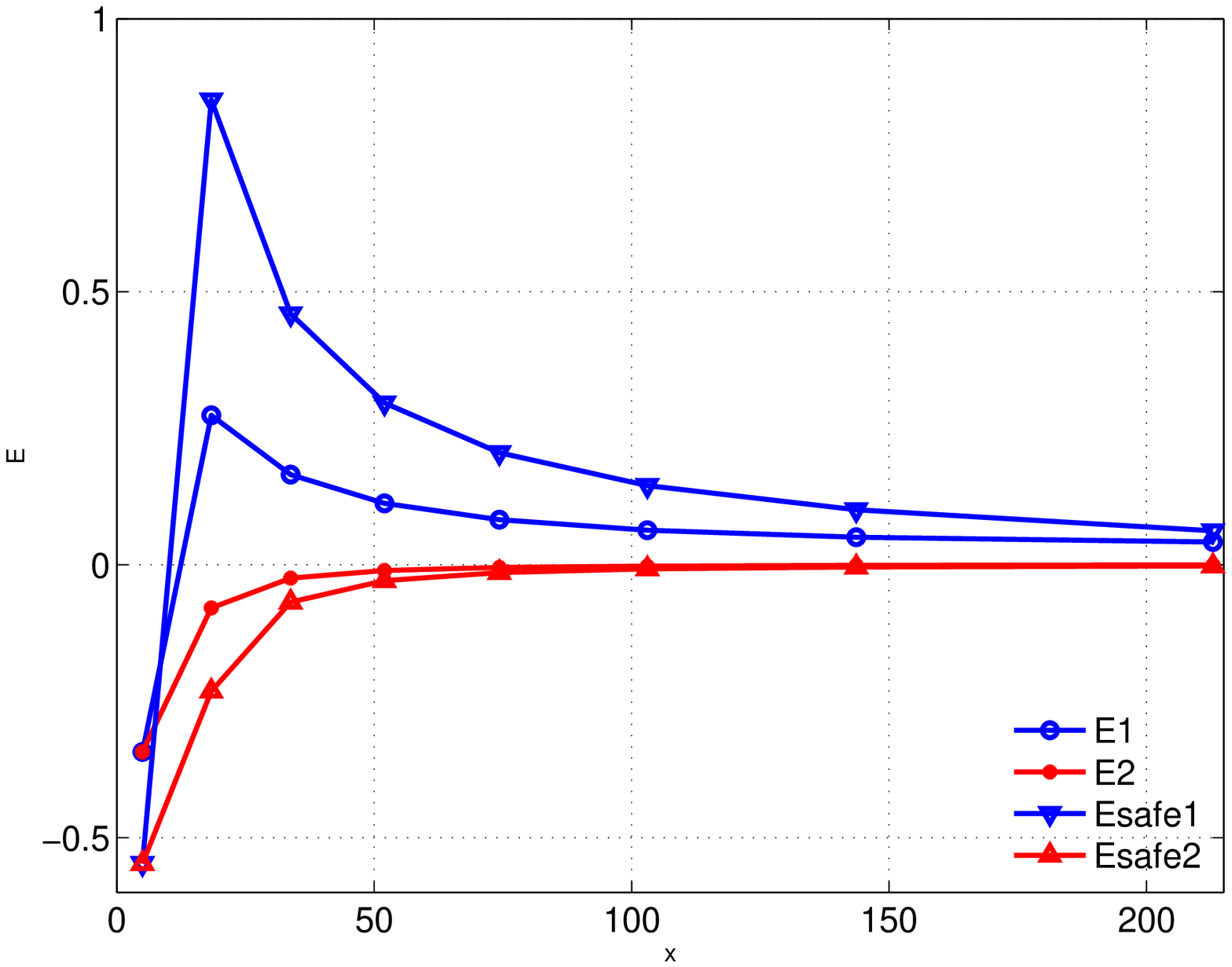} \\
\includegraphics[width=.7\textwidth]{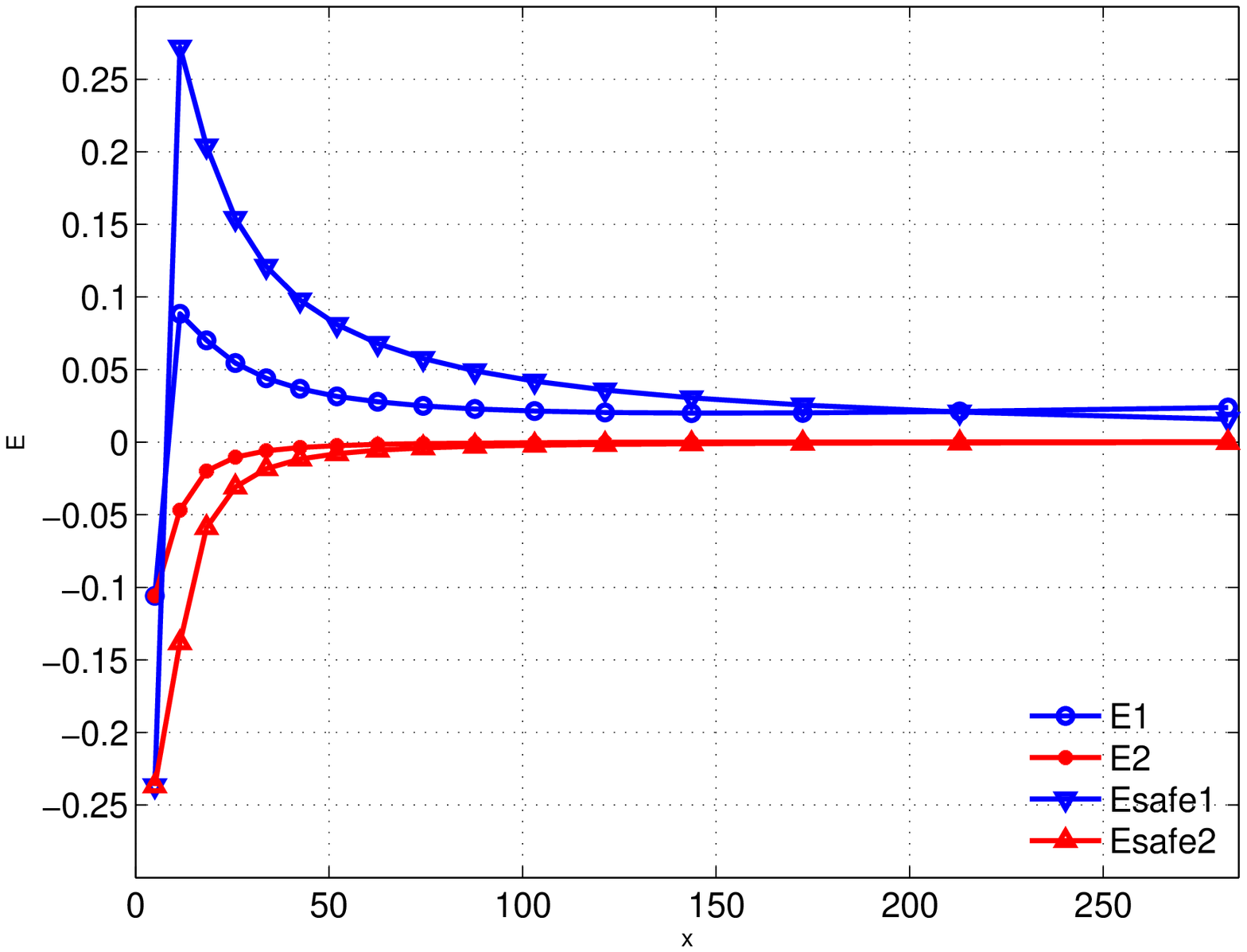} 
\caption{\it Safe error estimates provided by equation (\ref{eq:est2}) and true error by equation (\ref{eq:GE}) for the test problem (\ref{PPO:model}).
Here $E_1$ is the solution true error, $E_2$ true error for the solution first derivative, $Esafe_1$ the safe error estimate for the solution, and $Esafe_2$ the safe error estimate for the solution first derivative.}
\label{fig:PPO:error}
\end{figure}

Figure \ref{fig:PPO:error} shows two sample error estimates made by the error estimator (\ref{eq:est2}), from left to right we used $N=16$ and $N=32$.
It is easily seen that the safe estimator defined by equation (\ref{eq:est2}) provides upper bounds for the true error.

%In table \ref{tab:PPO} we report results for our test problem (\ref{PPO:model}) for different values of $N$.
%\input{TAB1.tex}
In table \ref{tab:PPO} we report, for different values of $N$, a few extrapolations for the free boundary value $R$ of our test problem (\ref{PPO:model}).
These values where computed according to the extrapolation formula
\begin{equation}\label{eq:Rextra:R}
R_{2N,k+1} = R_{2N,k} + \frac{R_{2N,k}-R_{N,k}}{2^{k+1}-1} \ , \quad \mbox{for} \ k = 0, 1, 2 \ .
\end{equation}
\begin{table}[!hbt]
\begin{center}{\renewcommand\arraystretch{1.2}
\begin{tabular}{r|cccc}
\hline
{$N$} & {$R_0$} & {$R_1$} & {$R_2$} & {$R_3$} \\
\hline  
%$10$   & $$ & $$ & $$ & $$ \\
$2$   & $2.1860735$ & $$ & $$ & $$ \\
$4$   & $3.2077586$ & $3.548320$ & $$ & $$ \\
$8$   & $4.1100007$ & $4.410748$ & $4.533952$ & $$ \\
$16$  & $4.6572864$ & $4.839715$ & $4.900996$ & $4.925466$ \\
$32$  & $4.8940600$ & $4.972985$ & $4.992024$ & $4.998093$ \\
$64$  & $4.9714936$ & $4.997305$ & $5.000779$ & $5.001363$ \\
$128$ & $4.9928640$ & $4.999987$ & $5.000370$ & $5.000343$ \\
$256$ & $4.9983436$ & $5.000170$ & $5.000196$ & $5.000184$ \\
$512$ & $4.9997042$ & $5.000158$ & $5.000156$ & $5.000153$ \\
\hline
\end{tabular}
}
\end{center}
\caption{Richardson's extrapolations for the free boundary value of the test problem (\ref{PPO:model}) with parameters fixed in (\ref{eq:parameters}).
Note that for this problem the exact value is $R = 5$.}
\label{tab:PPO}
\end{table}
In this table, since the values of $N$ can be seen on the first column, we omitted the first subscript for the notation defined in equation (\ref{eq:Rextra}) and used in equation (\ref{eq:Rextra:R}) for the extrapolated values of $R$.

%In the following 
\section{Conclusions}
In this paper we have presented the results obtained by a MATLAB version of a non-standard finite difference scheme for the numerical solution of the so-called
perpetual American put option model of financial markets and its generalizations.
This model can be derived from the celebrated Black-Scholes model letting the time goes to infinity.
Even in the classical Black-Scholes model, there is no known formula for the price of an American put with a finite exercise time (there are formulae for prices of infinite exercise time American put, American call option and European put and call options). 
A variety of numerical methods and approximations for the American put option price have been developed over the years. An overview of the various methods can be
found for example in Barone-Adesi \cite{Baroni:2005:SAP}.
The problem considered here is a \FBP \ defined on a semi-infinite interval, so that it is a non-linear problem complicated by a boundary condition at infinity.
By using non-uniform maps, we have shown how it is possible to apply the boundary condition at infinity exactly in contrast with the definition of a truncated boundary that introduces an error related to the replacement of infinity by a finite value, see for instance Nielsen et al. \cite{Nielsen:2002:PFF}.

As future work it would be relevant to extend our non-standard difference scheme to Black-Scholes models that are governed by partial differential equations defined on infinite domains.
Of course, we can apply the Landau's transform to the original moving boundary problem to get a problem defined on a fixed domain. 
The semi-discretization in time of the transformed problem with standard schemes like the first order Euler or high order Runge-Kutta type will result in a sequence of problems in the class (\ref{S:model}).

\vspace{1.5cm}

\noindent {\bf Acknowledgement.} {The research of this work was 
supported, in part, by the University of Messina and by the GNCS of INDAM.}


\begin{thebibliography}{10}

\bibitem{Amin:1994:CAO}
K.~Amin and A.~Khanna.
\newblock Convergence of {A}merican option values from discrete to
  continuous-time financial models.
\newblock {\em Math. Fin.}, 4:289--304, 1994.

\bibitem{Avellaneda:1996:MVR}
M.~Avellaneda and A.~Par\'as.
\newblock Managing the volatility risk of portfolios of derivative securities:
  the {L}agrangian uncertain volatility model.
\newblock {\em Appl. Math. Fin.}, 3:21--52, 1996.

\bibitem{Baroni:2005:SAP}
G.~Barone-Adesi.
\newblock The saga of the {A}merican put.
\newblock {\em J. Banking Finance}, 29:2909--2918, 2005.

\bibitem{Barraquand:1994:PAP}
J.~Barraquand and T.~Pudet.
\newblock Pricing of {A}merican path-dependent contingent claims, 1994.
\newblock Digital Research Laboratory, Paris.

\bibitem{Bensoussan:1984:TOP}
A.~Bensoussan.
\newblock On the theory of option pricing.
\newblock {\em Acta Appl. Math.}, 2:139--158, 1984.

\bibitem{Botta:1972:NSN}
E.~F.~F. Botta, D.~Dijkstra, and A.~E.~P. Veldman.
\newblock The numerical solution of the {N}avier-{S}tokes for laminar,
  incompressible flow past a parabolic cylinder.
\newblock {\em J. Eng. Math.}, 6:63--81, 1972.

\bibitem{Boyd:2001:CFS}
J.~P. Boyd.
\newblock {\em Chebyshev and Fourier Spectral Methods}.
\newblock Dover, New York, 2001.

\bibitem{Broadie:1996:MVR}
M.~Broadie and Detemple.
\newblock American option valuation: new bounds, approximations, and a
  comparison of existing methods.
\newblock {\em Rev. Financial Stud.}, 9:1211--1250, 1996.

\bibitem{Crank:FMB:1984}
J.~Crank.
\newblock {\em Free and Moving Boundary Problems}.
\newblock Clarendon Press, Oxford, 1984.

\bibitem{Davis:1972:NSN}
R.~T. Davis.
\newblock Numerical solution of the {N}avier-{S}tokes equations for symmetric
  laminar incompressible flow past a parabola.
\newblock {\em J. Fluid Mech.}, 51:417--433, 1972.

\bibitem{During:2012:HOC}
B.~D\"{u}ring and M.~Fourni\'e.
\newblock High-order compact finite difference scheme for option pricing in
  stochastic volatility models.
\newblock {\em J. Comput. Appl. Math.}, 236:4462--4473, 2012.

\bibitem{Elliot:1999:MFM}
R.~J. Elliot and P.~E. Kopp.
\newblock {\em Mathematics of Financial Markets}.
\newblock Springer, New York, 1999.

\bibitem{Fazio:2014:FDS}
R.~Fazio and A.~Jannelli.
\newblock Finite difference schemes on quasi-uniform grids for {BVP}s on
  infinite intervals.
\newblock {\em J. Comput. Appl. Math.}, 269:14--23, 2014.

\bibitem{Frey:2002:RMD}
R.~Frey and P.~Patie.
\newblock Risk management for derivatives in illiquid markets: a
  simulation-study.
\newblock In {\em Advances in Finance and Stoch.}, pages 137--159. Springer,
  Berlin, 2002.

\bibitem{Grosch:NSP:1977}
C.~E. Grosch and S.~A. Orszag.
\newblock Numerical solution of problems in unbounded regions: Coordinate
  transforms.
\newblock {\em J. Comput. Phys.}, 25:273--296, 1977.

\bibitem{Jandacka:2005:RAP}
M.~Janda\v{c}ka and D.~\v{S}ev\v{c}ovi\v{c}.
\newblock On the risk-adjusted pricing-methodology-based valuation of vanilla
  options and explanation of the volatility smile.
\newblock {\em J. Appl. Math}, 3:235--258, 2005.

\bibitem{Koleva:NSH:2006}
M.~N. Koleva.
\newblock Numerical solution of the heat equation in unbounded domains using
  quasi-uniform grids.
\newblock In I.~Lirkov, S.~Margenov, and J.~Wasniewski, editors, {\em
  Large-scale Scientific Computing}, volume 3743 of {\em Lecture Notes in
  Comput. Sci.}, pages 509--517, 2006.

\bibitem{Leland:1985:OPR}
H.~E. Leland.
\newblock Option pricing and replication with transactions costs.
\newblock {\em J. Finance}, 40:1283--1301, 1985.

\bibitem{Nilev:2013:EIS}
M.~Milev and A.~Tagliani.
\newblock Efficient implicit scheme with positivity preserving and smoothing
  properties.
\newblock {\em J. Comput. Appl. Math.}, 243:1--9, 2013.

\bibitem{Nielsen:2002:PFF}
B.~F. Nielsen, O.~Skavhaug, and A.~Tveito.
\newblock Penalty and front-fixing methods for the numerical solution of
  {A}merican option problems.
\newblock {\em J. Comput. Finance}, 5:69--97, 2002.

\bibitem{Richardson:1910:DAL}
L.~F. Richardson.
\newblock The approximate arithmetical solution by finite differences of
  physical problems involving differential equations, with an application to
  the stresses in a masonry dam.
\newblock {\em Proc. R. Soc. London Ser. A}, 210:307--357, 1910.

\bibitem{Richardson:1927:DAL}
L.~F. Richardson and J.~A. Gaunt.
\newblock The deferred approach to the limit.
\newblock {\em Proc. R. Soc. London Ser. A}, 226:299--349, 1927.

\bibitem{Schneider:GRR:2005}
F.~A. Schneider and C.~H. Marchi.
\newblock On the grid refinement ratio for ome-dimensional advection problems
  with nonuniform grids.
\newblock In {\em Proceedings of {COBEM} 2005}, 2005.
\newblock 18th International Congress of Mechamical Engineering, 8 pages.

\bibitem{vandeVooren:1970:NSS}
A.I. van~de Vooren and D.~Dijkstra.
\newblock The {N}avier-{S}tokes solution for laminar flow past a semi-infinite
  flat plate.
\newblock {\em J. Eng. Math.}, 4:9--27, 1970.

\bibitem{Veldam:PNG:1992}
A.~E.~P. Veldam and K.~Rinzema.
\newblock Playing with nonuniform grids.
\newblock {\em J. Eng. Math.}, 26:119--130, 1992.

\end{thebibliography}
\end{document}